\documentclass[pre,
twocolumn,
superscriptaddress,
]{revtex4-1}

\usepackage{times}
\usepackage{subcaption}
\usepackage{ragged2e}
\DeclareCaptionJustification{justified}{\justifying}
\usepackage{algorithm}
\usepackage{algpseudocode}
\usepackage{amsmath}
\usepackage{amssymb}
\usepackage{amsthm}
\usepackage{bbm}
\usepackage{enumerate}
\usepackage{enumitem}
\usepackage{pgfplots, pgfplotstable}
\usepackage{placeins}
\usepackage[titletoc,title]{appendix}
\usepackage{wrapfig,lipsum}
\usepackage{graphicx}
\usepackage{dcolumn}
\usepackage{bm}
\usepackage{overpic}
\usepackage{multirow}
\usepackage{tabularx,booktabs}
\usepackage{url}

\usepackage{float}
\floatstyle{plaintop}
\restylefloat{table}

\algrenewcommand\algorithmicindent{1.0em}%
\makeatletter
\newcommand{\algmargin}{\the\ALG@thistlm}   
\makeatother
\algnewcommand{\parState}[1]{\State%
\parbox[t]{\dimexpr\linewidth-\algmargin}{\strut #1\strut}}

\begin{document}

\title{A feasibility pump algorithm embedded in an annealing framework}
\author{Nicolas Pradignac, Maliheh Aramon} 
\affiliation{1QB Information Technologies Inc. (1QBit), 458-550 Burrard Street, Vancouver, BC V6C 2B5, Canada}

\author{Helmut G. Katzgraber}
\affiliation{Microsoft Quantum, Microsoft, Redmond, WA 98052, USA }
\affiliation{Department of Physics and Astronomy, Texas A\&M University, College Station, TX 77843-4242, USA}
\affiliation{Santa Fe Institute, 1399 Hyde Park Road, Santa Fe, NM 87501, USA}


\begin{abstract}
The feasibility pump algorithm is an efficient primal heuristic for finding feasible solutions 
to mixed-integer programming problems. The algorithm suffers mainly from fast convergence to local optima. 
In this paper, we investigate the effect of an alternative approach to circumvent this challenge by designing 
a two-stage approach that embeds the feasibility pump heuristic into an annealing framework. The 
algorithm dynamically breaks the discrete decision variables into two subsets based on the fractionality 
information obtained from prior runs, and enforces integrality on each subset separately. The feasibility pump algorithm 
iterates between rounding a fractional solution to one that is integral and projecting an infeasible integral solution 
onto the solution space of the relaxed mixed-integer programming problem. These two components are used in a Monte Carlo 
search framework to initially promote diversification and focus on intensification later. The computational results obtained from solving
$91$ mixed-binary problems demonstrate the superiority of our new approach over Feasibility Pump 2.0.
\end{abstract}
 
\date{\today}

\maketitle

\section{Introduction}
\label{sec:introduction}

Mixed-integer programming (MIP) is the most common approach for solving a wide range of optimization 
problems across a variety of domains, such as manufacturing~\cite{Pochet06}, transportation~\cite{Aramon13}, 
health care~\cite{Chan16}, energy~\cite{kuo16}, and finance~\cite{Mansini15}. Finding a feasible solution to an 
MIP problem is thus an important NP-complete problem and, in the last few decades, there have been many 
heuristic algorithms developed in the interest of more quickly finding a feasible solution to this problem. 
Examples of such heuristics are pivot-based~\cite{Balas80}, local search-based~\cite{Walser99}, and interior 
point-based~\cite{Orchers91} techniques. 

One of the MIP heuristic techniques embedded in most commercial state-of-the-art 
MIP solvers is the feasibility pump (FP) algorithm, a projection-based technique introduced by 
Fischetti et al.~\cite{Fischetti05}. The FP algorithm is an iterative algorithm that generates two solution 
sequences, one satisfying the linear constraints and the other satisfying the integrality 
constraints. The FP algorithm decrease the distance between the solutions generated in consecutive 
iterations until convergence. Despite there being no guarantee that a feasible solution will be obtained, 
it has been shown in practice that FP is very efficient \cite{Berthold17} and, as a result, many researchers 
have devised new methods for improving its performance.
 
The literature on the FP algorithm has focused on improving one of its three main components: rounding a 
fractional solution to an integral solution, projecting an integral solution onto the solution space subject to 
linear constraints, and perturbing the infeasible integral solution when the algorithm stalls.

In the original feasibility pump algorithm, a simple rounding procedure is used, where the decision variables 
are deterministically rounded to the nearest integer. Replacing this simple procedure with a randomized 
rounding strategy where each variable is randomly rounded to the nearest or second-nearest integer is 
among the early works on improving the rounding step of the algorithm~\cite{Bertacco07}. 
Bertacco et al.~\cite{Bertacco07} have further utilized their proposed rounding strategy in 
a two-stage approach for solving MIP problems with general-integer variables. The first stage 
concentrates on binary variables, relaxing the integrality of the other integer variables. The second 
stage restores the integrality of all integer variables and then initializes from the best solution 
found in the previous stage. 

To make the rounding step of the algorithm more efficient, Fischetti and Salvagnin~\cite{Fischetti09} 
have incorporated a constraints propagation mechanism~\cite{Savelsbergh94} to tighten the variable 
domains in order to generate an integral solution that is closer to feasibility. Their modified algorithm, 
called the Feasibility Pump~2.0 (FP2.0), is the most popular FP implementation. Their experimental results 
show that, not only does FP2.0 find a greater number of higher-quality feasible solutions than its predecessors, 
but it decreases the required computational time due to there being a reduced number of iterations. Examples of 
subsequent modifications to the rounding step are the work of Baena and Castro~\cite{Baena11} 
and Boland et al.~\cite{Boland14}, where the line search method is used to explore a large number of integer points. 
Variants of this approach have been proposed. It has been shown that the best variant, implemented on top of FP2.0, 
provides better-quality solutions than FP2.0 alone, but at the price of greater computational cost.

To enhance the algorithm's projection step, Achterberg and Berthold~\cite{Achterberg07} have 
modified the objective function of the projection model to be a convex combination of the distance 
function and the original MIP objective function. The convex coefficient is decreased at each iteration to 
gradually shift the focus from quality to integrality. This approach has been shown to improve the solution quality, 
but with a slight increase in the computational time. In a more recent work, Boland et al.~\cite{Boland12} 
have incorporated cutting planes in the FP algorithm to shrink the constraints polyhedron by removing the 
feasible solutions that are not integral. As the cutting plane mechanism could be expensive for very large 
instances of the MIP problem, the cuts are added only when the algorithm cycles, rather than at each iteration. 
It is also shown that the addition of cuts to FP2.0 boosts the solution quality; however, the process of generating cuts 
requires substantial computational effort, leading to an increase in computation time.

The infeasibility measure in the projection model of the original FP and FP2.0 algorithms is the $L_1$-norm 
between the rounded and projected solutions. De~Santis et al.~\cite{DeSantis10} have proposed 
replacing the common infeasibility measure by a weighted $L_1$-norm, in which the weights are calculated 
based on a nonlinear merit function. The motivation behind this modification is to consider the information 
from the projected solutions generated in previous runs and, consequently, to decrease the number of 
iterations needed to find a feasible solution. This version of the feasibility pump algorithm compares 
favourably with the original FP algorithm.

The main drawback of the FP algorithm is cycling, where the same sequence of points is visited repeatedly. 
Perturbation is the common approach used to address this drawback. Several perturbation mechanisms 
have been incorporated in FP2.0~\cite{Fischetti05}. Most of these techniques select the variables to be altered 
based on a  measure of fractionality. Dey et al.~\cite{Dey18} have recently introduced an improvement to 
the perturbation process. The new perturbation method is still based on fractionality, but it also has recourse 
to the support of the constraints that are violated by the current integral solution. Incorporating this perturbation 
method into FP2.0 results in a slight improvement in both the number of problems solved and the solution quality. 
A new framework based on the FP algorithm coupled to a biased random-key genetic algorithm has been 
recently developed by Andrade et al.~\cite{Andrade17}. The main idea behind the algorithm is to generate 
a pool of rounded and projected solutions, and merge them based on information, such as fractional infeasible 
solutions and their corresponding rounded solutions, collected from previous iterations. Their experimental 
results show that the proposed framework is efficient in solving hard MIP problems. It is able to find a greater 
number of feasible solutions than either FP2.0 or CPLEX. 

The literature on feasibility pumps has recently been expanded to include convex mixed-integer nonlinear 
programming problems~\cite{DAmbrosio10,DAmbrosio12}. We refer the interested reader to the 
survey by Berthold at al.~\cite{Berthold17}, which gives a comprehensive review of the literature on 
FP since its inception. 

As mentioned, the algorithms based on FP use perturbation to diversify the search space when 
the algorithm has been detected to be cycling. Such an approach is \emph{reactive} to escaping 
from local minima. Conversely, Monte Carlo simulation techniques such as simulated annealing \cite{Kirkpatrick83}, replica 
exchange \cite{Hukushima96,Predescu05}, and population annealing \cite{HuIb03,Ma10,WaMaKa15b} utilize randomization \emph{proactively} through an aggressive 
diversification performed at the beginning of the search, then gradually guiding the search toward the 
region that has a higher potential of containing the global optimum. In this paper, we propose a new 
two-stage approach that embeds the feasibility pump algorithm in a simulated annealing 
framework~\cite{Kirkpatrick83}, called the Annealed Feasibility Pump (AFP). The new approach 
divides the variables dynamically into two subsets based on their fractionality values, which have been 
collected from previous runs. It then enforces integrality on each subset separately, where the sequences 
of rounding and projecting are utilized in an annealing framework to initially promote diversification, and focus 
on intensification later. 

The paper is organized as follows. We present an overview of the FP algorithm in Section~\ref{sec:fp}. 
Then, in Section~\ref{sec:fpsa}, we describe the AFP algorithm. Computational results comparing the 
performance of the new algorithm with the state-of-the-art feasibility pump algorithms are presented 
and discussed in Section~\ref{sec:computational_results}. We present our conclusion in Section~\ref{sec:conclusion}.

\section{The Feasibility Pump Algorithm}
\label{sec:fp}

In this section, we provide an overview of the FP algorithm and present a modification of the projection step.

\subsection{Overview}
\label{sec:overview_fp}

The feasibility pump algorithm finds a feasible solution to a mixed-integer programming (MIP) problem, 
defined as 
\begin{align} 
\centering
\mbox{(MIP)~~~~~~~~} \min \{c^Tx| Ax \leq b\text{, }x_i \in \mathbb{Z}~ ~\forall i \in I\}, & \label{eq:MIP}
\end{align}
\noindent where $c \in \mathbb{R}^n$, $A \in \mathbb{R}^{m \times n}$, $b \in \mathbb{R}^m$, and 
$I \subseteq \{1, 2, \ldots, n\}$. The set of variable indices restricted to being integral is denoted by $I$. The upper and 
lower bounds on variables are included in the constraint matrix $A$. Let $\mathcal{P} = \{x: Ax \leq b\}$ 
be the polyhedron that includes all the LP-feasible solutions in which the linear constraints are satisfied 
and the integrality constraints on the variables in set $I$ are relaxed~\cite{Fischetti05,Bertacco07,Fischetti09,Boland14,Achterberg07,Boland12}. 
In case all the integral decision variables are binary, we refer to the problem as a mixed-binary programming (MBP) problem.

The FP algorithm is an iterative approach that alternates between sequences of LP-feasible solutions 
that are not integral and LP-infeasible solutions that are integral. In this paper, we refer to the former 
as ``relaxed'' solutions and the latter as ``integral'' solutions. The goal of the algorithm is to 
decrease the distance between solutions in consecutive iterations until converging to an LP-feasible, 
integral solution. More specifically, the algorithm starts from the optimal solution to the relaxed MIP model, 
where the integrality constraints are relaxed and the linear constraints are satisfied. The LP-feasible 
solution is then rounded, forming an integral solution that might violate some of the linear constraints. 
To finish one iteration of the algorithm, the rounded solution is projected back onto the space 
$\mathcal{P}$ to reach the next LP-feasible solution. The rounding and projecting procedure 
is repeated until a feasible solution to the MIP problem formulated as expression~\eqref{eq:MIP} is found. 
The main drawback of this algorithm is that it might stall, that is, the current integral solution has 
been visited in a previous iteration. Perturbation is the common approach incorporated into FP, forcing 
the algorithm to explore a new region of the search space. Algorithm~\ref{alg:fp} summarizes the main 
steps of the FP algorithm, where $\bar{x}$ and $\tilde{x}$ represent the current relaxed solution and 
the current integral solution, respectively.

\begin{algorithm}[H]
\caption{\footnotesize{\,\,The FP Algorithm}}
\begin{algorithmic}[1]
\For {each run}
    \State {solve the relaxed MIP problem to obtain $\bar{x}$}
    \State {round $\bar{x}$ to obtain $\tilde{x}$}
    \While {stopping conditions not reached}
        \State {project $\tilde{x}$ back onto $\mathcal{P}$ to obtain the next $\bar{x}$}
        \State {round $\bar{x}$ to obtain the next $\tilde{x}$}
        \If {$\tilde{x}$ has been already visited}
            \If{cycle length is 1}
                \State {perform a weak perturbation on $\tilde{x}$}
            \Else 
                \State {perform a strong perturbation on $\tilde{x}$}
            \EndIf
        \EndIf
    \EndWhile
\EndFor
\end{algorithmic}
\label{alg:fp}
\end{algorithm}

The  key components of the FP algorithm are as follows.

\vspace{0.3cm}
\paragraph*{Rounding} To transform the relaxed solution into an integral one, rounding $\bar{x}_i$, the $i$-th element 
of $\bar{x}$, to the nearest integer is the simplest procedure found in the literature~\cite{Fischetti05}. 
In this paper, we apply the randomized rounding procedure where each variable 
of the relaxed solution is randomly rounded to the nearest or second-nearest integer~\cite{Bertacco07}. 
The details of this procedure are explained in Section~\ref{sec:fpsa_}. 

\vspace{0.3cm}
\paragraph*{Projection} This step finds a new point in the space $\mathcal{P}$, minimizing the objective 
feasibility pump function~\cite{Achterberg07}:
\begin{align}
\centering
~\min_{x} ~ & (1-\alpha) \frac{\Delta(x, \tilde{x})}{\left\lVert \Delta \right\Vert} + \alpha \frac{c^Tx}{\left\lVert c \right\Vert}  \label{feasibility_pump_obj} \\
\mathrm{s.t.} ~& ~~ Ax \leq b & \notag
\end{align}

Here, $\Delta(x, \tilde{x}) = \sum_{i \in I} |x_i - \tilde{x}_i|$ is the $L_1$-norm distance between 
a point \mbox{$x \in \mathcal{P}$} and a given integral point $\tilde{x}$. In our implementation, 
we have converted the distance function $\Delta(x, \tilde{x})$ into a linear expression by adding auxiliary 
variables and extra constraints~\cite{Fischetti05}. Further, $\left\lVert \Delta \right\Vert = \sqrt{|I|}$, and $\left\lVert c \right\Vert$ is the $L_2$-norm 
of vector $c$. The coefficient $\alpha$ is updated at the end of each iteration, 
adjusting the contribution of quality versus integrality. Initially, it is set at one to promote 
high-quality solutions and is decreased at each iteration in order to favour integrality. 

\begin{figure}
\centering
\begin{subfigure}[b]{0.5\textwidth}
\centering
\includegraphics[width=\textwidth]{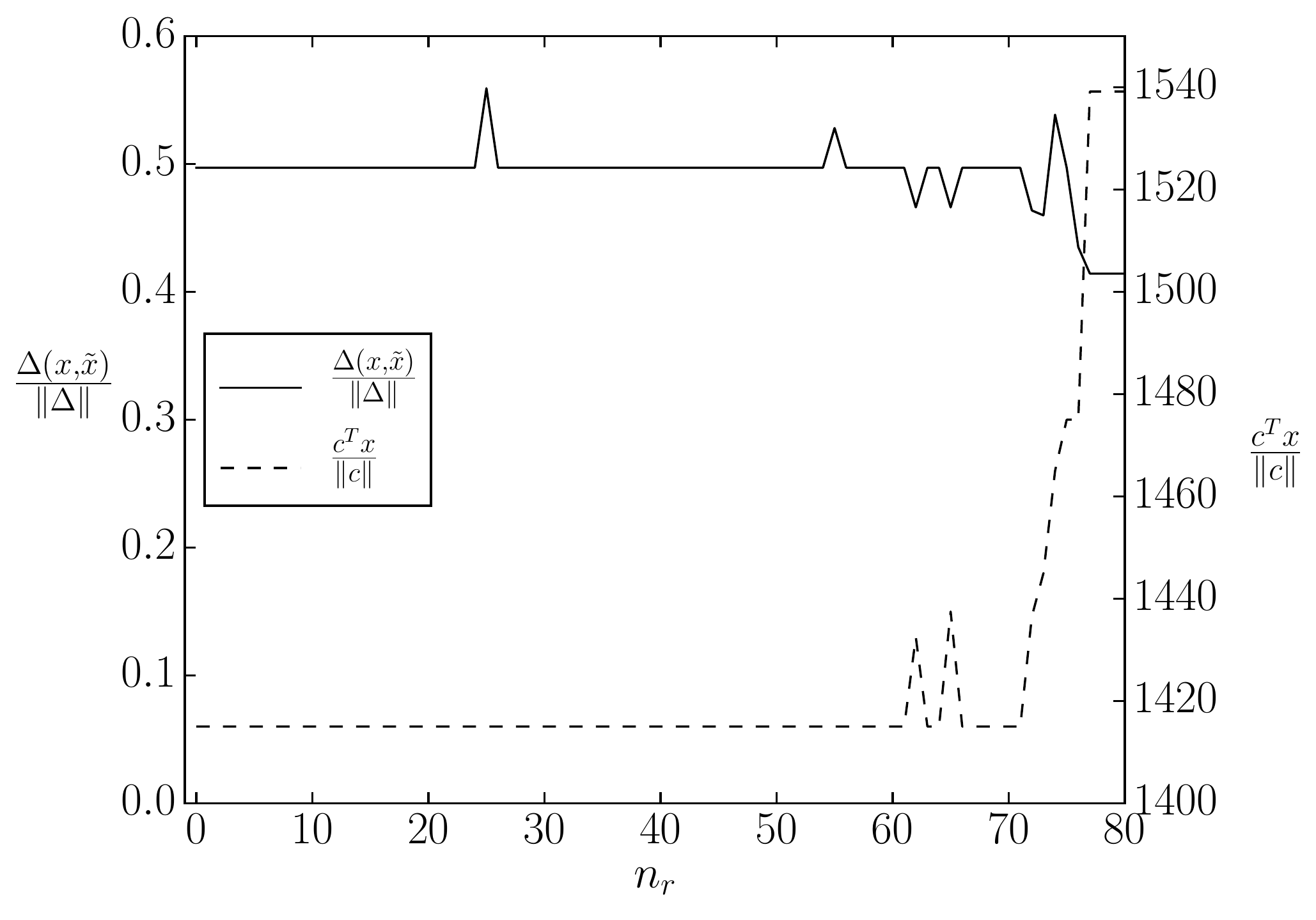}
\caption{}
\label{fig:misc07_distance_term_norm}
\end{subfigure}
\begin{subfigure}[b]{0.5\textwidth}  
\centering 
\includegraphics[width=\textwidth]{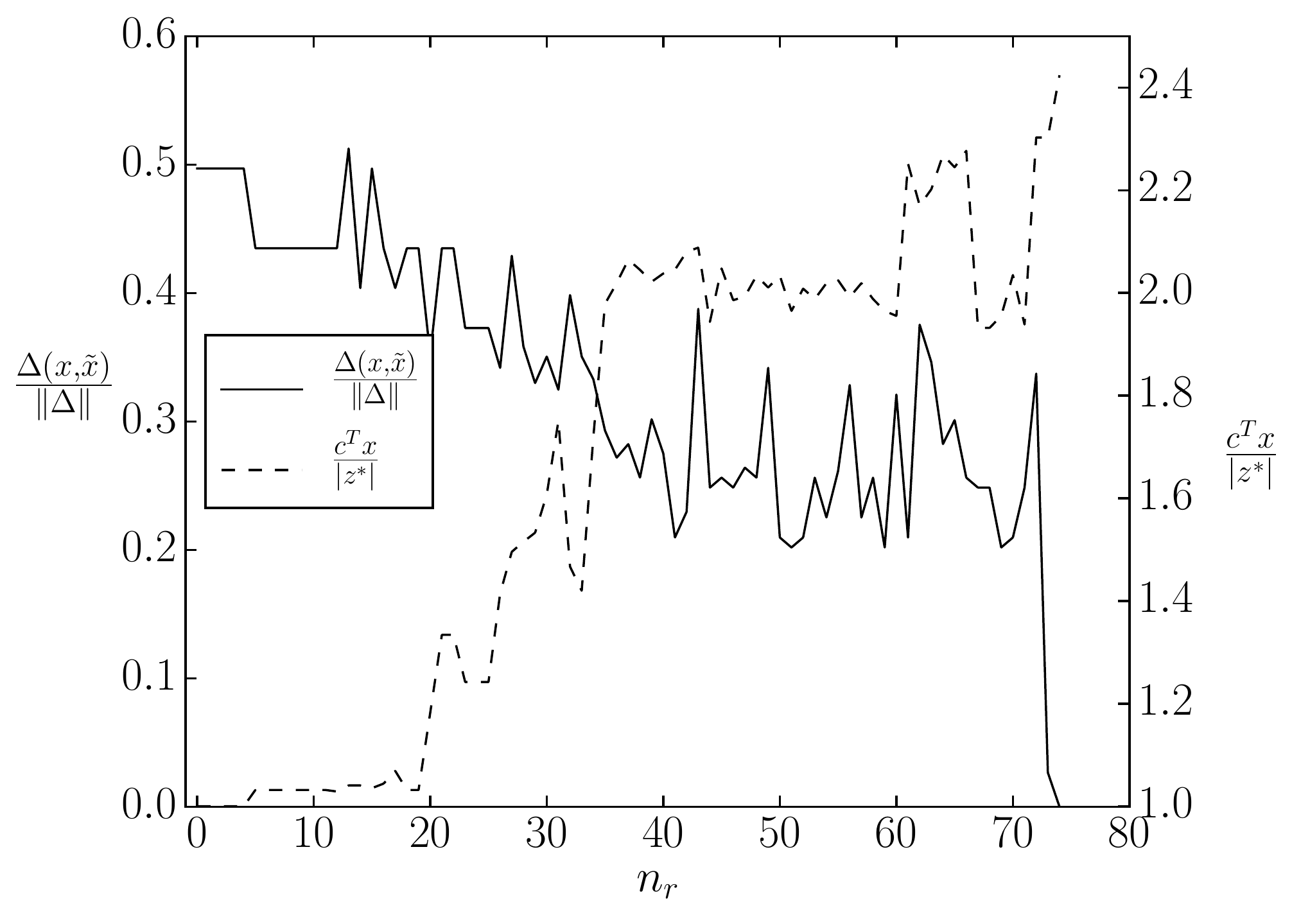}
\caption{}
 \label{fig:misc07_quality_term_optimal_relaxed}
\end{subfigure}
 \caption{The distance and quality values through iterations ($n_r$) of an example run of the benchmarking problem 
 \texttt{misc07} are shown in (a) and (b) where, in respective terms, the quality $c^Tx$ is normalized by 
 the coefficient vector norm $\left\lVert c \right\Vert$ and by the optimal objective function of the relaxed MIP model $z^*$.}
\label{fig:misc07_quality_distance_term}
\end{figure}

The main idea in dividing the terms of expression~\eqref{feasibility_pump_obj} 
by their Euclidean norms is to ensure that both the distance and quality terms are non-dominant 
in the projection step. However, in our experimentation on problem \texttt{misc07} from the \texttt{MIP2003} 
library \cite{MIPLIB10}, we discovered that normalizing the quality 
term $c^Tx$ by the coefficient vector norm does not guarantee a balanced contribution 
from both quality and distance terms. As shown in Fig.~\ref{fig:misc07_distance_term_norm}, the distance values 
on the left-hand side vertical axis are almost constant across the iterations in the range $(0.4, 0.5)$. 
However, the quality values, shown on the right-hand side vertical axis, are on the order of thousands, resulting in an 
approximately three-orders-of-magnitude difference in the range of distance and quality 
values. In this example, the relaxed solution $\bar{x}$ remains predominantly unchanged across the iterations. 
Although the same integral solution was visited repeatedly by the algorithm, no cycle was detected and the perturbation 
procedure was not performed because the decrease in $\alpha$ was greater than $\delta = 0.005$ in two consecutive runs~\cite{Achterberg07}. Such a phenomenon 
usually occurs when the objective function of the MIP problem, that is, $c^Tx$ contains a few unbounded continuous 
decision variables and $\left\lVert c \right\Vert$ is very small. To address this issue,  Achterberg et al.~\cite{Achterberg07} 
have used the MIP preprocessing techniques embedded in CPLEX prior to running the FP algorithm. In this paper, 
we take a different approach. Specifically, we divide the quality term by the optimal objective function value 
of the relaxed MIP model $z^*$ in order to balance the contributions of both terms. Fig.~\ref{fig:misc07_quality_term_optimal_relaxed}
 illustrates the effect of the new normalization factor where the distance term has a non-monotonic decreasing trend 
 until it reaches zero and the distance and quality values are more balanced through the iterations. 

\paragraph*{Perturbation} At the end of each iteration $t$, the algorithm checks whether the current 
integral solution $\tilde{x}$ has been visited in previous iterations. If it has been visited at iteration $t-1$, a cycle of 
length one is detected and the variables of $\tilde{x}$ undergo a weak perturbation~\cite{Fischetti05}. Otherwise, a strong 
perturbation is performed during which a larger number of variables are altered~\cite{Fischetti05}. 

\section{The Annealed Feasibility Pump Algorithm}
\label{sec:fpsa}

Simulated annealing (SA)~\cite{Kirkpatrick83} is an iterative, generic heuristic method that 
emulates the process of physical annealing to find the global optimum of a combinatorial 
optimization problem. Considering a finite solution space, SA is a discrete time Markov chain 
that acts on a nonincreasing temperature schedule. At each temperature value, the probability 
of transitioning to a neighbouring state is calculated according to a Metropolis criterion which 
promotes diversification at high temperatures and intensification at low temperatures. As with the 
majority of heuristic algorithms, reaching a local optimum at the lowest temperature is a likely scenario in SA. 
To increase the likelihood of reaching the global minimum, this process is repeated multiple times, each time 
starting from a random initial solution. In this paper, we present an algorithm we have designed called 
the Annealed Feasibility Pump (AFP), which utilizes the main components of the feasibility pump algorithm, 
including rounding and projection, in a simulated annealing framework.

In what follows, we first describe the AFP algorithm and then present several modifications 
we have incorporated into the algorithm to further improve its performance. 

\subsection{The Annealed Feasibility Pump Algorithm}
\label{sec:fpsa_}
The AFP algorithm embeds key features of the feasibility pump algorithm into a simulated annealing framework. 
The key steps of AFP are given in Algorithm~\ref{alg:fpsa}.

\begin{algorithm}[H]
\caption{\footnotesize{\,\,The AFP Algorithm}}
\begin{algorithmic}[1]
\For {each run}
    \State {initialize to a random initial pair ($\bar{x}, \tilde{x}$)}
    \While {stopping conditions not reached}
    	\For {each $\alpha$ and neighbourhood function}
        		\State {generate a new solution pair and compute its \\ 
		           ~~~~~~~~~~~~~fractionality}
        		\State {compute $\Delta d$, the decrease in the fractionality}
		\If{$\Delta d \leq 0$}
			\State {update ($\bar{x}, \tilde{x}$)} 
            		\State {break}
		\ElsIf {$\mbox{U}(0, 1) \leq e^{\frac{- \Delta d}{\alpha}}$}
            		\State {update ($\bar{x}, \tilde{x}$)} 
		        \State {break}
        		\EndIf
    	\EndFor
    \State {update $\alpha$}
\EndWhile
\EndFor
\end{algorithmic}
\label{alg:fpsa}
\end{algorithm}

The algorithm iteratively generates new solution pairs~($\bar{x}$, $\tilde{x}$). Each pair, in respective terms, 
is composed of a relaxed solution and an integral solution. The algorithm starts from an initial pair of a 
random integral solution and its corresponding projected relaxed solution. At each consecutive iteration, 
a new solution pair is generated by one of the neighbourhood functions. The algorithm moves to the new solution pair according 
to a Metropolis criterion. This process is repeated until a stopping condition has been reached.

\vspace{0.3cm}
\paragraph*{Initial Pair} The algorithm starts from an initial integral
solution $\tilde{x}$, where a random integer value is assigned to each integer variable 
within its bounds, and the continuous variables are set to their corresponding values of the relaxed model's optimal solution. 
Projecting the initial random solution $\tilde{x}$  back onto the LP-feasible solution space $\mathcal{P}$ 
results in an LP-feasible solution $\bar{x}$, generating an  initial solution pair. 

\vspace{0.3cm}
\paragraph*{Neighbourhood Functions} For each value of $\alpha$, the algorithm randomly iterates over a list of  
neighbourhood functions until a new solution pair is accepted or the list is exhausted. The list includes 
five  neighbourhoods that differ mainly in terms of the number of variables they alter. Each neighbourhood 
function generates a new solution pair by first generating a new integral solution based on the current solution pair 
and then projecting it back onto the space $\mathcal{P}$ to obtain a new relaxed solution. 

\begin{enumerate}
	\item \emph{Randomized Rounding}: This function aims to generate a new integral solution close 
	to the current relaxed solution~\cite{Bertacco07}. For each discrete variable~$i$, a value  $\tau(\omega)$ 
	is computed as shown, where $\omega$ is a random variable from the uniform distribution $\rm{U}(0,1)$: 
	\begin{align}
		\tau(\omega) = 
			\begin{cases}
				 2\omega(1-\omega),& \text{if } \omega \leq 0.5 \\
				1 - 2\omega(1-\omega),& \text{if } \omega > 0.5 \notag
			\end{cases}
	\end{align}
	The value of the variable $i$ in the new integral solution, $\tilde{x}^{\mbox{\scriptsize{n}}}_i$, is then equal to 
	$\lfloor \bar{x}^{\mbox{\scriptsize{c}}}_i + \tau(\omega) \rfloor$. The random variable $\tau(\omega)$ is bounded 
	between $[0, 1]$ with a mode of $0.5$. When $\tau(\omega) = 0.5$, the variable $i$ is rounded to the nearest integer.
	\\
	\item \emph{Weak Perturbation}: This function is utilized when the original MIP problem is an MBP 
	problem. Given the ordered list of binary variables with positive fractionality 
	in the current relaxed solution $\bar{x}^{\mbox{\scriptsize{c}}}$, a subset of variables with the highest 
	fractionality values is randomly selected to be flipped~\cite{Fischetti05}. The number of variables flipped is 
	generated from a discrete uniform distribution ${\mbox{U}}(T/2 , 3T/2)$, where $T$ is equal to a fraction 
	of the length of the ordered list defined earlier. We have used $0.1$ as the fraction value in our experimentation. 
	\\
	\item \emph{Strong Perturbation}: This function is utilized only for the MBP problems where 
	a random variable $\omega_i$ is first generated from ${\mbox{U}}(l,u)$ for each binary variable $i$. 
	If the sum of the fractionality of the variable $i$ in the current relaxed solution and $\max(0, \omega)$ 
	is larger than $0.5$, then the variable is flipped~\cite{Fischetti05}. We have used ${\mbox{U}}(-0.3,0.7)$ in our experiments.
	\\
	\item \emph{Weak Perturbation Domain}: This function is similar to the weak perturbation neighbourhood and is 
	applied to both MIP and MBP problems. The variables to be perturbed are selected following similar steps, 
	with the difference that the fractionality value for the discrete variable $i$ is normalized by its 
	domain size $D_i$. Defining $w_i$ as the window size for each variable $i$, the new integral solution is generated according to
	\begin{align}
		\tilde{x}^{\mbox{\scriptsize{n}}}_i = 
			\begin{cases}
				{\mbox{U}}(\bar{x}^{\mbox{\scriptsize{c}}}_i, \bar{x}^{\mbox{\scriptsize{c}}}_i + w_i), 
				& \text{if } \bar{x}^{\mbox{\scriptsize{c}}}_i \geq  \tilde{x}^{\mbox{\scriptsize{c}}}_i  \\
				{\mbox{U}}(\bar{x}^{\mbox{\scriptsize{c}}}_i - w_i, \bar{x}^{\mbox{\scriptsize{c}}}_i), 
				& \text{otherwise }  \notag
			\end{cases}
	\end{align}

	We have used $w_i = \max(50, 0.05\times D_i)$ in our experiments. We have further ensured that 
	the new chosen values respect variables' upper and lower bounds. 
	\\
	\item \emph{Strong Perturbation Domain}: This function, motivated by the work of Fischetti and Salvagnin~\cite{Fischetti09}, 
	is applied to both MIP and MBP problems where the values of half of the discrete decision variables chosen 
	at random are perturbed. Defining $D_i$ and $w_i$, respectively, as the variable $i$'s domain size and the window size, and 
	further letting $U_i$ and $L_i$ be the variable $i$'s upper and lower bounds, respectively, the new value for $i$ is
	\begin{align}
		\tilde{x}^{\mbox{\scriptsize{n}}}_i = 
			\begin{cases}
				{\mbox{U}}(\bar{x}^{\mbox{\scriptsize{c}}}_i, U_i), 
				& \text{if } D_i < a \text{ and } \bar{x}^{\rm{\scriptsize{c}}}_i \geq  \tilde{x}^{\mbox{\scriptsize{c}}}_i  \\
				{\mbox{U}}(L_i, \bar{x}^{\mbox{\scriptsize{c}}}_i), &  \text{if } D_i < a \text{ and } \bar{x}^{\mbox{\scriptsize{c}}}_i < \tilde{x}^{\mbox{\scriptsize{c}}}_i  \\
				{\mbox{U}}(U_i - w_i, U_i), & \text{if } U_i - \tilde{x}^{\mbox{\scriptsize{c}}}_i \leq \beta \times D_i \\
				{\mbox{U}}(L_i, L_i + w_i), &  \text{if } \tilde{x}^{\mbox{\scriptsize{c}}}_i - L_i \leq \beta \times D_i \\
				{\mbox{U}}(\bar{x}^{\mbox{\scriptsize{c}}}_i - w_i, \bar{x}^{\mbox{\scriptsize{c}}}_i + w_i), & \text{otherwise } \notag
			\end{cases}
	\end{align}
	We have set $a = 10$ and $\beta = 0.1$ in our experiments. In each of the above cases, we ensure that the new selected value 
	is within the variable's bounds.	 
\end{enumerate}

\vspace{0.3cm}
\paragraph*{Metropolis Criterion} Once a new solution pair has been generated by a neighbourhood function, 
the $L_1$-norm distance between $\bar{x}^{\mbox{\scriptsize{n}}}$ and $\tilde{x}^{\mbox{\scriptsize{n}}}$ 
of the new pair, that is, the fractionality, is calculated. We further define $\Delta d$ to 
be the difference between the fractionality of the new solution pair and the current pair. 
The algorithm then proceeds from the current pair $(\bar{x}^{\mbox{\scriptsize{c}}}, \tilde{x}^{\mbox{\scriptsize{c}}})$ 
to the new solution pair $(\bar{x}^{\mbox{\scriptsize{n}}}, \tilde{x}^{\mbox{\scriptsize{n}}})$ 
according to the following probability: 
\begin{align}
\mathbbm{P}\big((\bar{x}^{\mbox{\scriptsize{c}}}, \tilde{x}^{\mbox{\scriptsize{c}}}) 
\rightarrow (\bar{x}^{\mbox{\scriptsize{n}}}, \tilde{x}^{\mbox{\scriptsize{n}}}) \big) = 
\begin{cases}
        1,& \text{if } \Delta d \leq 0 \\
        0,& \text{if } \alpha = 0 \\
        e^{\frac{- \Delta d}{\alpha}} ,& \text{otherwise} \notag
\end{cases} 
\end{align}
\noindent{If $\Delta d \leq 0$, the new solution pair is  accepted as it is closer to being an integral, LP-feasible 
solution pair. Otherwise, the move is accepted with a probability that decreases as the search 
progresses.}

\vspace{0.3cm}
\paragraph*{Stopping Conditions}
Each run of the algorithm terminates as soon as one of the following conditions is satisfied: reaching 
an integral, LP-feasible solution; exceeding the maximum number of iterations per run $n_r$; or 
not improving the fractionality measure in the last $k$ iterations. At the end of each run, the 
acceptance probability is updated before initiating the next run (see Section \ref{subsec:fpsa_modification} for additional details). 
The entire process terminates when the total number of iterations, $n_t$, has been exhausted. 

\subsection{Modifications}
\label{subsec:fpsa_modification}

In this section, we present several modifications incorporated into the AFP algorithm 
to further improve its performance either in terms of the number of feasible solutions found or the
solution quality.

\subsubsection{Dynamic Acceptance Probability}
\label{subsubsec:dynamic_accep_prob}

The coefficient $\alpha$ in the projection model plays the role of temperature  
in the AFP algorithm. Under a proper temperature schedule in a simulated annealing framework, 
the probability of accepting a move that worsens the solution is high at high temperatures 
and decreases, reaching zero, as the temperature is lowered. To guarantee such 
behaviour in the AFP algorithm and, taking into account that the $\alpha$ coefficient 
is bounded in $[0,1]$, the $\Delta d$ values need to be normalized by the normalization 
factor $\left\lVert \Delta d \right\Vert$. Note that, in the AFP algorithm, a move is considered  
a ``bad move'' if it increases the current fractionality value, yielding $\Delta d > 0$. 

The normalization factor is computed based on the data collected in the first run of the AFP algorithm. 
Initially, $\left\lVert \Delta d \right\Vert$ is set to one, and the worst move with the largest $\Delta d_{\rm{max}}$ is 
tracked. At the end of the first run, $\left\lVert \Delta d \right\Vert$ is updated such that the probability 
of accepting the worst move at a predetermined high $\alpha_{\rm{h}}$ is equal to $p_{\rm{h}}$. More 
specifically, we have found $\left\lVert \Delta d \right\Vert$ such that

\begin{align}
e^{\frac{- \Delta d_{\rm{max}}}{\left\lVert \Delta d \right\Vert \alpha_{\rm{h}}}} = p_{\rm{h}}. & \notag 
\end{align}

\begin{figure}
	\begin{subfigure}[b]{0.5\textwidth} 
		\centering
		\includegraphics[width=\textwidth]{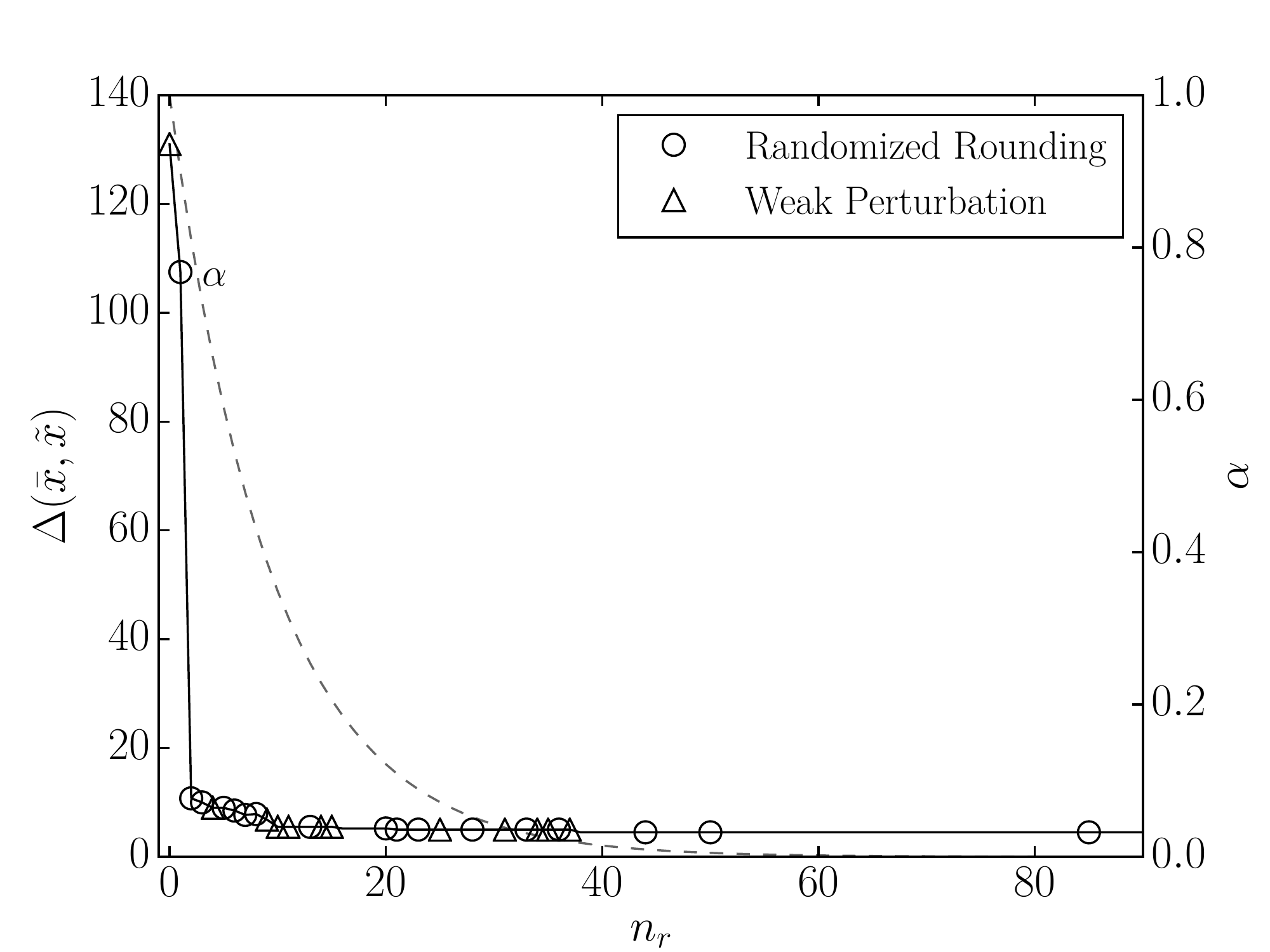}
		\caption{}
		\label{fig:misc_fractionality_without_normalization}
	\end{subfigure}
         \begin{subfigure}[b]{0.5\textwidth}  
        		 \centering
		\includegraphics[width=\textwidth]{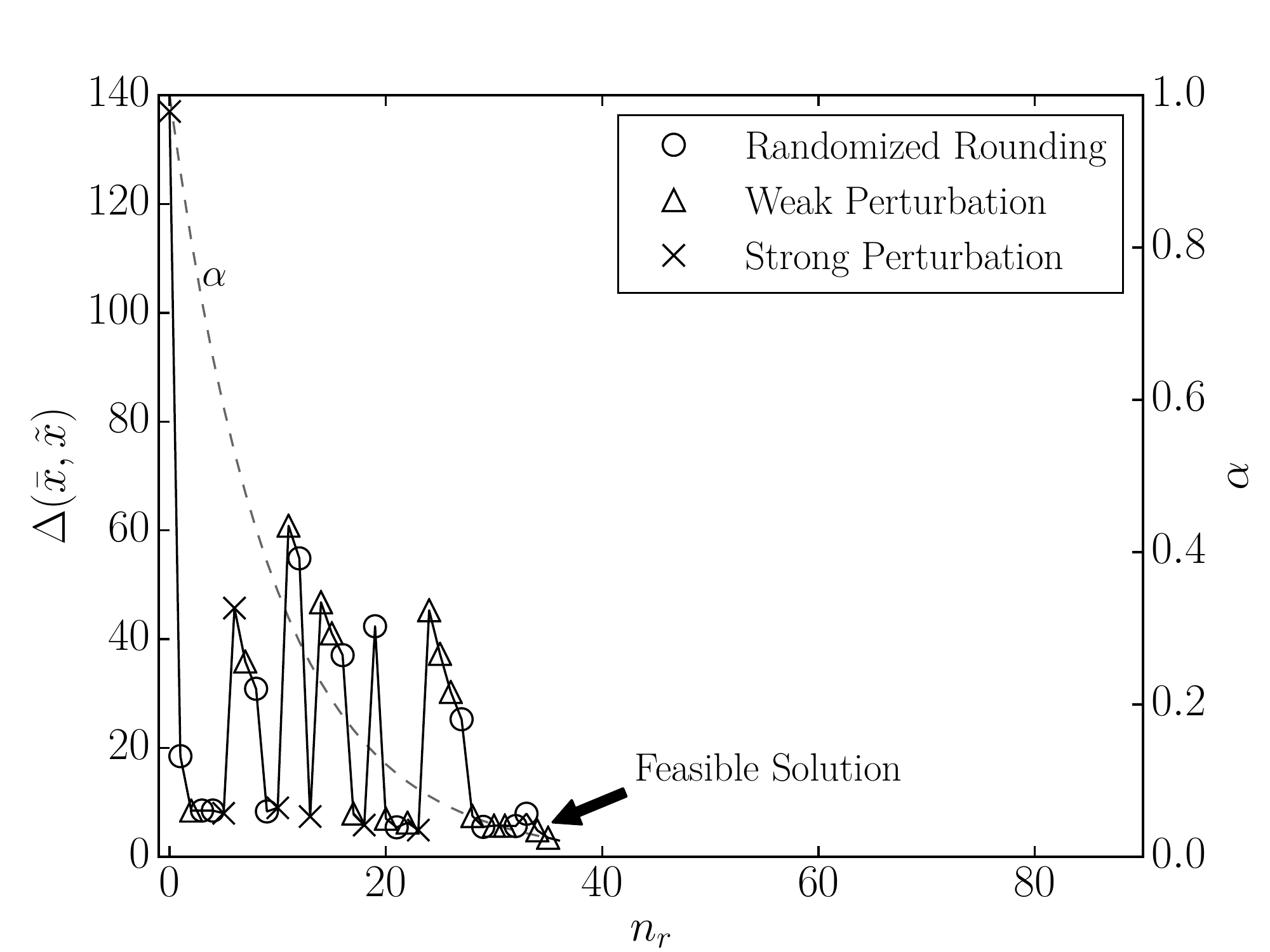}
		\caption{}
		\label{fig:misc_fractionality}
	\end{subfigure}
\caption{Fractionality value of the current solution ($\Delta(\bar{x}, \tilde{x})$) obtained by the AFP algorithm 
at each iteration ($n_r$) in an example run of benchmarking problem \texttt{misc07}. The $y$-axis on the right-hand sides 
of the graphs represents the corresponding $\alpha$ values, and the marker at each iteration 
indicates the accepted neighbourhood function. The acceptance probabilities have been calculated 
(a) without normalizing the fractionality values, and (b) with the fractionalilty values normalized.}	
\label{fig:misc}
\end{figure}

Fig.~\ref{fig:misc} illustrates how the fractionality value of the current solution of the AFP algorithm 
changes over iterations in an example run of the benchmarking problem \texttt{misc07}. 
If the fractionality values are not normalized, the probability of accepting bad moves is very low even 
at high $\alpha$ values, and the simulated annealing framework essentially acts as a greedy heuristic 
(see Fig.~\ref{fig:misc_fractionality_without_normalization}). Normalizing the fractionality values 
obviates such behaviour, as shown in Fig.~\ref{fig:misc_fractionality}. The probability 
of accepting bad moves is high at early iterations and, as a result, the fractionality value oscillates, 
reaching different regions of the solution space. As the algorithm progresses, the fractionality value stabilizes, 
intensifying the search to improve upon the current solution. Furthermore, the example run depicted 
in Fig.~\ref{fig:misc_fractionality_without_normalization}, terminates because the fractionality measure 
has not improved in the last 70 iterations, whereas in Fig.~\ref{fig:misc_fractionality}, the run ends, 
as a feasible solution has been found.

At the end of each run, the normalization factor is recalculated based on whether a feasible solution 
has been found. If a feasible solution has not been found, the acceptance probability $p_{\rm{h}}$ is 
increased to boost diversification in the subsequent run. As soon as a feasible solution is found 
during a run, the probability $p_{\rm{h}}$ is decreased to encourage intensification, increasing the 
odds of reaching a higher-quality feasible solution.

\subsubsection{Hard-Variable Fixing}
\label{sec:hard_variable_fixining}

Motivated by the common rule of branching on those variables that have the smallest domain first, 
Bertacco et al. \cite{Bertacco07} developed a two-stage approach for solving MIP problems with 
general integer variables using the FP algorithm as explained in Section~\ref{sec:introduction}. 
The main idea behind their two-stage approach is that enforcing integrality on the binary 
variables first would likely result in the integrality of the majority of the general integer variables. 
The second-stage problem would then be the easier problem of deriving the integrality of the remaining integer variables. 

In this paper, we present a new two-stage approach that differs in terms of the variables 
included in the first stage. Let us define $\mathcal{B}$ and $\mathcal {I}$ as the set of binary 
and general-integer variables, respectively. Further, for each discrete decision variable $i \in \mathcal{B} \cup \mathcal {I}$, 
we define $r_i$ as the infeasibility rank that counts the number of times the variable $i$ is fractional at the end of a run. 
We refer to the variables with the largest $r_i$ values as ``hard variables'' and concentrate on enforcing 
integrality on these variables in the first stage. 

\begin{algorithm}[H]
\caption{\footnotesize{\,\,The AFP Algorithm with Hard-Variable Fixing}}
\begin{algorithmic}[1]
\State {$\mathcal{H} = \emptyset, n_h=0, r_i = 0 ~(\forall i \in \mathcal{B} \cup \mathcal {I})$}
\If{$n_h = 0$} \label{initial_if}
	\State{run AFP algorithm} \label{initial_fpsa_run}
	\If{$\bar{x}$ is feasible}
		\State{return to step \ref{initial_fpsa_run}}
	\Else
		\State{run Algorithm \ref{alg:HardVariables}}
		\State{go to step \ref{stage1}}
	\EndIf
\Else
	\State{stage 1: relax easy variables \label{stage1}}
	\State{stage 1: run AFP algorithm}
	\If{$\bar{x}$ is feasible}
		\State{stage 2: fix variables in set $\mathcal{H}$ at values found}
		\State{stage 2: run AFP algorithm}
		\State{run Algorithm \ref{alg:HardVariables}}
		\State{return to step \ref{stage1}}
	\Else
		\State{run Algorithm \ref{alg:HardVariables}}
		\State{return to step \ref{stage1}}
	\EndIf
\EndIf
\end{algorithmic}
\label{alg:hard_variable_fixing}
\end{algorithm}

Our algorithm (see Algorithm~\ref{alg:hard_variable_fixing}) first solves the MIP problem as defined 
in expression~\eqref{eq:MIP}, where all discrete decision variables are included in the AFP algorithm. After 
reaching the stopping condition for a run, the infeasibility rank value $r_i$ is calculated for 
each discrete decision variable $i \in \mathcal{B} \cup \mathcal {I}$ if a feasible solution is 
not found. We further define $\mathcal{H}$ as the hard-variables set that includes $n_h$ discrete decision variables with the largest 
infeasibility rank. We refer to the remaining decision variables as ``easy variables''. 
Note that, at the end of the first run, the infeasibility rank of the decision variables is either one or zero, and 
the set $\mathcal{H}$ might contain a mix of binary and integer decision variables.

\begin{algorithm}[H]
\caption{\footnotesize{\,\,Updating the Hard-Variables Set $\mathcal{H}$}}
\begin{algorithmic}[1]
\If{$n_h = 0$}
	\State {run Algorithm \ref{alg:infeasibilityRanking}}
	\State{$\mathcal{H} = \{i| r_i > 0\}$}
	\State{$n_h = |\mathcal{H}|$}
\Else
	\If{Algorithm \ref{alg:hard_variable_fixing} is at the end of stage $1$ and $\bar{x}$ is \\
	   \hspace{0.75cm} infeasible}
		\State{$n_h = 0.8 \times n_h$}
		\State{run Algorithm \ref{alg:infeasibilityRanking}}
		\State{$\mathcal{H} \leftarrow$ the $n_h$ variables with the highest infeasibility \\ 
		\hspace{1.6cm} ranks}
	\ElsIf{Algorithm \ref{alg:hard_variable_fixing} is at the end of stage $2$ and $\bar{x}$ is \\ 
		\hspace{1.3cm} infeasible}
		\State{$n_h = 1.2 \times n_h$}
		\State{run Algorithm \ref{alg:infeasibilityRanking}}
		\State{$\mathcal{H} \leftarrow$ the $n_h$ variables with the highest infeasibility \\
		  \hspace{1.6cm} ranks}
	\ElsIf{Algorithm \ref{alg:hard_variable_fixing} is at the end of stage $2$ and $\bar{x}$ is \\
	         \hspace{1.3cm} feasible}
		\State{$n_h = 0.8 \times n_h$}
		\State{run Algorithm \ref{alg:infeasibilityRanking}}
		\State{$\mathcal{H} \leftarrow$ the $n_h$ variables with the highest infeasibility \\
		            \hspace{1.6cm} ranks}
	\EndIf	
\EndIf
\end{algorithmic}
\label{alg:HardVariables}
\end{algorithm}

\begin{algorithm}[H]
\caption{\footnotesize{\,\,Updating the Infeasibility Rank $r$ }}
\begin{algorithmic}[1]
\State {infeasible variables $\leftarrow$ integer variables with a positive 
	  fractionality in $\bar{x}$}
\For {each infeasible variable $i$}
	\State {$r_i \mathrel{+}= 1$} 
\EndFor
\end{algorithmic}
\label{alg:infeasibilityRanking}
\end{algorithm}

Having defined the set $\mathcal{H}$, the algorithm begins the first stage, where the integrality constraints 
on the easy variables are relaxed and the AFP algorithm aims to find an LP-feasible solution 
that is integral with respect to the variables in set $\mathcal{H}$. After the first stage, if the 
algorithm has successfully reached its aim, it proceeds to the second stage. Otherwise, the 
infeasibility rank of the variables in $\mathcal{H}$ is updated and the  set $\mathcal{H}$ 
is redefined by decreasing $n_h$. The algorithm then attempts to solve a new modified problem in the first stage 
with fewer discrete decision variables.

 In the second stage, the decision variables in $\mathcal{H}$ are fixed at the values found in the 
 first stage and the  goal is to enforce integrality on the easy variables. At the end of this stage, 
 if a feasible solution has been found, it means that a feasible solution to the original MIP problem 
 in expression~\eqref{eq:MIP} has been found. The algorithm then repeats the entire process, where 
 the set $\mathcal{H}$ is updated by decreasing the number of hard variables $n_h$. The main purpose 
 of adjusting the set $\mathcal{H}$ is to focus on finding a better-quality solution in the remaining runs. 
 In case the second stage does not reach integrality on all easy variables, the algorithm restarts 
 at the first stage, where the number of hard variables $n_h$ is increased in the hope that 
 the new first-stage problem will guide the algorithm in the direction of obtaining integrality. 

The stopping conditions in each stage are similar to the ones defined in Section~\ref{sec:fpsa_} for each 
run of the AFP algorithm, and the algorithm terminates when the total number of iterations $n_t$ is reached. 

\section{Computational Results}
\label{sec:computational_results}

In this section, we discuss the performance of the AFP algorithm. First, we investigate 
the effect of embedding the simulated annealing framework in the feasibility pump 
for deriving diversification and intensification. We then compare the AFP algorithm to the current state-of-the-art 
feasibility pump solvers such as FP2.0~\cite{Fischetti09} and FPC~\cite{Boland12}. Finally, we evaluate the impact 
of a two-stage approach that uses hard-variable fixing.  

Our testbed consists of a subset of MBP instances used in~\cite{Boland12} from 
the \texttt{COR@L} \cite{CORAL} and \texttt{MIP2003} \cite{MIPLIB03} libraries. The instances 
are classified according to three classes, A, B, and C, in increasing order of difficulty based 
on the number of major restarts required to obtain an integer, LP-feasible solution~\cite{Boland12}. 
That is, class C contains the hardest instances which have a ragged solution space filled with many local minima. 
The AFP and FP algorithms have been implemented in Python, with Gurobi 8.1 used as the linear programming solver. The 
number of problem instances in classes A, B, and C are $31$, $30$, and $30$, respectively.  
The performance metrics used to compare algorithms are the following:

\begin{enumerate}
	\item [] {\bf{gap}} The percentage difference between the best integral, LP-feasible solution found $(\bar{x}^*)$ 
				   and the optimal or the best known solution $(x^*)$:
				   \begin{align}
				   	100 \times \frac{c^T\bar{x}^* - c^Tx^*}{\max(|c^Tx^*|, 1)} & \notag
				 \end{align}
	\item [] {\bf{success}} The ratio between the number of successful runs and the total number of runs. A run 
					 is considered successful if it finds a feasible integral solution. 
	\item [] {\bf{time}} The time taken per run.
\end{enumerate}

\noindent {\bf{AFP vs. FP }} The AFP and FP algorithms have been run in two stages with the hard-variable fixing 
modification, where the projection model's objective function is normalized by the optimal objective function of 
the relaxed MIP model. In both algorithms, the same stopping conditions are used: the total number 
of iterations $n_t$, set at $5000$; the maximum number of iterations per run $n_r$, set at $150$; 
and the number of iterations $k$ to improve the fractionality, set at $70$. The initial parameter $\alpha$ 
in both algorithms is set to one and decreases geometrically with a factor of $0.9$ per iteration. The 
list of neighbourhood functions in the AFP algorithm includes randomized rounding, 
weak perturbation domain, and strong perturbation domain. We experimented with different combinations of 
neighbourhood functions on a small subset of problem instances and found that the combination 
of these three neighbourhood functions resulted in better performance. The parameter $p_h$ used in the 
Metropolis criterion is initially set at $0.7$. It decreases or increases, in respective terms, 
to $0.9 \times p_h$ or $\sqrt{p_h}$ at each run depending on whether a feasible solution has been found in the previous run. 

Table \ref{table:fpsa_fp} reports, for each algorithm, the number of problem instances with lower gaps 
($\mathcal{N}^{\tiny{\mbox{gap}}}$), the average success rates, and the average time per run. 
Detailed results for each problem are provided in Tables~\ref{table:classA},~\ref{table:classB},~and~\ref{table:classC}. 
We know from the literature that the MIP preprocessing techniques available in the state-of-the-art MIP solvers are used prior 
to running feasibility pump algorithms and can reduce the total running time of FP-type algorithms~\cite{Achterberg07}. 
However, we do not apply Gurobi's presolving techniques before running the AFP and the FP algorithms here.

The AFP algorithm outperforms FP in its yielding of higher-quality solutions. Both algorithms achieve 
similar solution quality on $11$, $7$, and $8$ problem instances of  classes A, B, and C, respectively. 
The performance difference between AFP and FP is the highest for class A and the lowest for class C, 
as expected. Since the key difference between AFP and FP is AFP's utilization of  simulated annealing, 
the results in Table \ref{table:fpsa_fp} indicate that using a different strategy for diversifying the search space 
is more effective in reaching high-quality solutions than finding a feasible solution for hard problems. 

\begin{table}
\caption{The number of problems with lower gaps ($\mathcal{N}^{\tiny{\mbox{gap}}}$), the average success rate, 
and the average time per run for the AFP and FP algorithms across problem classes. The AFP and 
 FP algorithms are run in two stages with hard-variable fixing.}
\label{table:fpsa_fp}
\centering
\begin{tabular}{ll|rr}  \hline \hline
\multicolumn{2}{c|}{Metrics} 		 			& AFP 		& FP 	 \\ \hline
\multirow{3}{*}{$\mathcal{N}^{\tiny{\mbox{gap}}}$} 	& Class A 		&     	\textbf{15}	&	 5		\\ 
									   	& Class B 		& 	\textbf{13}	&	10		\\ 
									   	& Class C 	&	\textbf{12}	&	10		\\ \hline
\multirow{3}{*}{success (avg.)} 					& Class A 		&     	0.79	&	0.79		\\ 
									   	& Class B 		& 	0.71	&	\textbf{0.77}		\\ 
									   	& Class C 	&	\textbf{0.49}	&    0.45		\\ \hline
\multirow{3}{*}{time (avg.)} 					& Class A 		&    \textbf{28.06} &    29.48			\\ 
									   	& Class B 		&    14.30	&    \textbf{6.96}		\\ 
									   	& Class C 	&  \textbf{267.04}	&  423.73		\\ \hline \hline
\end{tabular}			
\end{table}									   

\vspace{0.5cm}							   								   
{\bf{AFP vs. FP2.0 and FPC }} Table~\ref{table:fpsa_fp2_fpc} shows a comparison of the AFP algorithm with two 
state-of-the-art feasibility pump solvers,  FP2.0 and FPC, with respect to the number of instances solved to a 
higher-solution quality (i.e., lower gap). We rely on the results presented by Boland et al.~\cite{Boland12} for the state-of-the-art solvers 
(\mbox{Tables}~\ref{table:classA},~\ref{table:classB},~and~\ref{table:classC}) give detailed results for each problem).

As mentioned in Section~\ref{sec:introduction}, FP2.0 is an extension of the feasibility 
pump algorithm proposed by Fischetti and Salvagnin~\cite{Fischetti09} that replaces the 
simple rounding step with a strategy based on constraint propagation, where rounding 
is performed after inferring tighter bounds on each variable. More recently, 
Boland et al.~\cite{Boland12} designed a new algorithm called FPC that utilizes 
cutting planes in the perturbation step of FP2.0 to remove the LP-feasible, non-integral local minimum solutions.

\begin{table}
\caption{The number of problem instances solved to a better level of quality using each solver 
for each problem class. The comparison is performed pairwise between AFP and FP2.0, and 
between AFP and FPC.}
\label{table:fpsa_fp2_fpc}
\centering
\begin{tabular}{l|rrr|r} \hline \hline
 Algorithm 	& Class A 		& Class B 		& Class C 	& Total \\\hline 
AFP  		& 10 			& 17 			& 	15		& \textbf{42} \\
FP2.0 		& 11			&  9 			&  	11		& 31 \\ \hline \hline
AFP  		& 6 			& 11 			&  	10		& 27 \\
FPC 			& 15 			& 15 			&  	20		& \textbf{50} \\ \hline \hline
\end{tabular}
\end{table}

As shown in Table~\ref{table:fpsa_fp2_fpc}, AFP provides better-quality solutions than FP2.0 for a larger number of problem instances: $42$ instances versus $31$. 
The performance difference is more significant 
for harder problems, that is, those from class B and class C. However,  the FPC algorithm clearly outperforms 
AFP, as it achieves higher-quality solutions for two times as many problem instances. 
Our early experimentation on MIP problems with general-integer variables 
showed that FP2.0 and FPC perform better than the AFP algorithm.

Our current results demonstrate that using a Monte Carlo simulation framework for diversifying the search space achieves superiority 
over FP2.0 in solving MBP problems. In these binary problems, the constraint propagation is not very effective because the rounding concerns changing only the direction, not deciding on a new value. In this case, the FP2.0 
is essentially similar to the original FP algorithm~\cite{Fischetti05}. We believe that embedding the FPC algorithm in a simulated annealing 
framework, where the cutting planes are incorporated into a new neighbourhood function, has a high potential for 
achieving solutions of improved quality for both MBP and MIP problems. We plan to consider this approach in 
 future work.

The detailed results in Table~\ref{table:classA} further show that for problem instance \texttt{neos-934278}, 
AFP finds a solution of better quality than the best known solution, which can be obtained from the \texttt{COR@L} library. 

\vspace{0.5cm}
{\bf{Impact of Hard-Variable Fixing}} To understand the importance of the two-stage approach 
with hard-variable fixing, Table~\ref{table:fpsa_fp_one} shows a comparison of both AFP and FP utilized in a two-stage approach with 
a one-stage approach. The two-stage approach has a clear advantage over the one-stage 
approach for both AFP and FP. The performance difference becomes 
more pronounced in solving the problems of class C, the hardest problem class.

\begin{table}[h]
\caption{The number of problem instances solved to a higher level of quality using each solver 
in each problem class. The comparison is performed pairwise between the two-stage 
AFP algorithm and the one-stage AFP algorithm (AFP-1), and between the two-stage FP algorithm 
and the one-stage FP algorithm (FP-1).}
\label{table:fpsa_fp_one}
\centering
\begin{tabular}{l|rrr|r} \hline \hline
 Algorithm 	& Class A 		& Class B 		& Class C 	& Total \\\hline 
AFP  		&  	6		& 	13		& 	14		& \textbf{33} \\
AFP-1 		& 	6		&   	11		&  	6		&  22 \\ \hline \hline
FP  			&  	6		&  	11		&  	14		&  \textbf{30} \\
FP-1 		&  	6		&  	11		&  	8		&  25 \\\hline \hline
\end{tabular}
\end{table}

\section{Conclusion}
\label{sec:conclusion}

Perturbing the current integral solution is a common randomization approach  
to diversify search in the feasibility pump algorithm when it stalls. Such an 
approach  to the cycling challenge is reactive, addressing it after it occurs. In this paper, we have proposed a new two-stage 
approach that uses randomization proactively by embedding the feasibility pump algorithm in a Monte 
Carlo simulation framework.

The algorithm dynamically divides the discrete decision variables into two  subsets based on 
the information from prior runs, and enforces integrality separately on each subset. At each iteration, 
the algorithm performs a random walk, moving to a new solution pair composed of a relaxed solution and an integral 
solution based on a Metropolis criterion. The probability of moving to a pair with a high fractionality value is 
initially high, then decreases gradually as the search progresses toward finding a feasible, integral solution. 

Our computational results for $91$ mixed-binary optimization problems show the benefit of 
the new diversification scheme in reaching higher-quality solutions. Our results further indicate that  
the proposed approach is superior to FP2.0, especially in solving harder problem instances with ragged solution spaces 
filled with many local minima. 

In addition, our results show that the original feasibility pump algorithm embedded in the simulated annealing 
framework underperforms FP2.0 in solving mixed-integer optimization problems with general 
integer variables. Potential ideas for future work include coupling the FP2.0 algorithm with the simulated annealing framework, utilizing 
Monte Carlo replica exchange techniques for a better exploration of the search, and designing 
new neighbourhood functions that incorporate cutting planes in the projection step.

\begin{acknowledgements}
The authors would like to thank Jon Machta for useful discussions, Ehsan Iranmanesh for 
reviewing the manuscript, and Marko Bucyk for editorial help. Part of the work of H.~G.~K.'s research 
is based upon work supported in part by the Office of the Director of National Intelligence (ODNI),
Intelligence Advanced Research Projects Activity (IARPA), via MIT Lincoln Laboratory Air Force 
Contract No.~FA8721-05-C-0002. The views and conclusions contained herein are those of the 
authors and should not be interpreted as necessarily representing the official policies or endorsements, 
either expressed or implied, of ODNI, IARPA, or the U.S.~Government. The U.S.~Government is authorized 
to reproduce and distribute reprints for Governmental purpose notwithstanding any copyright annotation thereon.
\end{acknowledgements}

\bibliographystyle{apsrevtitle}
\bibliography{refs}

\begin{table*}[hb]
\caption{Detailed results of different solvers on class A problems. The ``inf'' means that the 
problem has not been solved to feasibility.}
\label{table:classA}
\centering
\begin{tabular}{|l|rrrrrr|rrrr|rrrr|}  \hline
\multirow{2}{*}{{\bf{Class A}}} 	& \multicolumn{6}{c|}{{\bf{gap}}} &  \multicolumn{4}{c|}{{\bf{success}}}	&  \multicolumn{4}{c|}{{\bf{time}}} 	 \\  
					& AFP & FP & AFP-1  & FP-1 & FP2.0 & FPC 
					& AFP & FP & AFP-1 & FP-1 & 
					AFP & FP & AFP-1 & FP-1 \\  \hline
\texttt{aflow40b} 		&    54.45	&  103.85	&  110.02	&  113.78	&    16.90	&    15.30	& 	0.33	&      0.71	& 	0.18	&      0.80	&     19.54	& 	7.77	&   23.60	& 	6.70	 \\ 
\texttt{bc1} 			& 	3.29	& 	3.59	& 	3.29	&      3.59	& 	3.10	& 	3.20	& 	0.99	&      1.00	& 	1.00	& 	1.00	&     10.68	& 	6.33	&   15.12	& 	5.78	 \\ 
\texttt{dano3\_3} 		& 	0.00	& 	0.00	& 	0.00	& 	0.00	& 	0.00	& 	0.00	& 	1.00	& 	1.00	& 	1.00	& 	1.00	& 	9.27	& 	9.00	&   16.36	& 	9.01	 \\
\texttt{fiber} 			&  213.22	&  226.71	&  188.31	&  226.71	&  133.10	& 	3.80	& 	0.99	&      1.00	& 	1.00	& 	1.00	& 	2.52	& 	1.21	&    2.56	& 	1.16	 \\
\texttt{leo1} 			&    22.96	&    23.20	&    21.01	&    25.33	&    17.60	&     11.90	& 	0.99	&      1.00	& 	1.00	& 	1.00	&     12.68	&     11.24	&   13.17	&     12.23	 \\
\texttt{mkc} 			&    35.75	&    46.15	&    35.70	&    46.15	&    49.80	&     13.10	& 	0.99	& 	1.00	& 	1.00	& 	1.00	& 	8.20	& 	7.02	&     8.62	& 	7.07	 \\
\texttt{neos-1200887}  	& 	2.70	& 	2.70	&      2.70	&      2.70	& 	5.40	& 	2.70	& 	1.00	& 	1.00	& 	1.00	& 	1.00	& 	0.13	& 	0.21	& 	0.14	& 	0.18	 \\
\texttt{neos-495307}    	&	0.06	& 	0.14	& 	0.06	& 	0.14	& 	5.70	& 	5.70	& 	1.00	& 	1.00	& 	1.00	& 	1.00	&    12.07	&     16.30	&     13.67	&     16.25	 \\
\texttt{neos-504815}    	&  272.07	&  340.65	&  224.04	&  375.40	&    42.00	&     10.70	& 	0.72	& 	1.00	& 	0.76	& 	1.00	& 	2.86	& 	0.50	&      2.48	& 	0.45	 \\
\texttt{neos-525149}    	&	0.00	& 	0.00	& 	0.00	& 	0.00	& 	0.00	& 	0.00	& 	1.00	& 	1.00	& 	1.00	& 	1.00	&     11.13	&    13.69	&    19.43	&     13.68	 \\
\texttt{neos-547911}    	&	0.00	& 	0.00	& 	0.00	& 	0.00	& 	7.70	& 	7.70	& 	1.00	& 	1.00	& 	0.99	& 	1.00	& 	6.10	& 	4.09	& 	6.84	& 	4.52	 \\
\texttt{neos-584146}    	&	0.00	& 	0.00	& 	0.00	& 	0.00	& 	0.00	& 	0.00	& 	1.00	& 	1.00	& 	1.00	& 	1.00	& 	0.52	& 	0.27	& 	0.50	& 	0.25	 \\
\texttt{neos-612125}    	&	0.00	& 	6.79	& 	0.12	& 	6.79	& 	0.50	& 	0.00	& 	1.00	& 	1.00	& 	1.00	& 	1.00	&  	2.29	& 	3.71	& 	4.72	& 	4.34	 \\ 
\texttt{neos-691058}    	& 	0.00	& 	0.00	& 	0.00	& 	0.00	& 	0.00	& 	0.30	& 	0.79	& 	0.87	& 	0.81	& 	0.81	&    13.03	& 	5.61	&     12.55	&      3.39	 \\
\texttt{neos-791021}    	&    40.00	&    53.33	& 	inf	&    40.00	&    40.00	& 	6.70	& 	0.43	& 	0.44	& 	0.00	& 	0.30	&     71.61	&     60.39	&     86.00	&     69.26	 \\
\texttt{neos-810286}    	& 	3.23	& 	inf	& 	0.00	& 	0.07	& 	3.90	& 	0.50	& 	0.25	& 	0.00	& 	0.37	& 	0.20	&     42.09	&     37.45	&     43.72	&    17.50 \\
\texttt{neos-825075}    	& 	0.00	& 	0.00	& 	0.00	& 	0.00	& 	0.00	& 	0.00	& 	0.94	& 	0.99	& 	0.96	& 	1.00	& 	0.96	& 	0.30	& 	0.73	& 	0.27	 \\
\texttt{neos-827175}    	&	6.67	& 	8.04	& 	6.67	& 	8.04	& 	0.00	& 	0.00	& 	1.00	& 	1.00	& 	1.00	& 	1.00	&     65.41	&  302.09	&     88.62	&   301.75	 \\
\texttt{neos-847302}    	& 	0.00	& 	0.00	& 	0.00	& 	0.00	& 	0.00	& 	0.00	& 	0.94	& 	0.97	& 	0.96	& 	0.94	& 	2.64	& 	0.98	& 	2.44	& 	0.94	 \\
\texttt{neos-906865}    	& 	0.00	& 	0.03	& 	0.00	& 	0.03	& 	0.00	& 	0.00	& 	1.00	& 	1.00	& 	1.00	& 	1.00	& 	0.82	& 	0.53	& 	0.75	& 	0.48	 \\
\texttt{neos-933815}    	& 	2.17	& 	2.87	& 	3.72	& 	2.74	& 	5.50	& 	4.10	& 	0.99	& 	1.00	& 	1.00	& 	1.00	& 	1.57	& 	1.00	& 	1.57	& 	0.96	 \\
\texttt{neos-934278}    	&     -1.08	&     -2.87	&     -1.80	&     -2.87	&    25.20	&     26.60	& 	1.00	& 	1.00	& 	1.00	& 	1.00	&    53.23	&    96.09	&   105.28	&   107.76	 \\
\texttt{neos-935674}    	&    50.00	&   -16.67	&    83.33	&    33.33	& 1.73e+3	&   850.00	& 	0.78	& 	0.73	& 	0.79	& 	0.78	& 	7.78	& 	6.27	& 	6.75	& 	3.68	 \\
\texttt{neos-941782}    	&    16.67	&     11.11	&    22.22	& 	5.56	& 	5.60	& 	0.00	& 	0.14	& 	0.20	& 	0.45	& 	0.67	&     5.36	& 	6.52	& 	8.79	& 	2.51	 \\
\texttt{neos-957323}    	& 	0.00	& 	0.00	& 	0.00	& 	0.00	&      0.00	& 	0.00	& 	1.00	& 	1.00	& 	1.00	& 	1.00	&   209.87	& 105.90	&  318.21	&   378.07 		 \\
\texttt{neos12} 			&  120.27	& 	inf	&  117.57	&  117.57	&      7.70	& 	7.70	& 	0.58	& 	0.00	& 	0.65	& 	0.63	&    49.59	&     72.07	&    46.04	&     20.05	 \\
\texttt{neos5}			& 	0.00	& 	0.00	& 	0.00	& 	0.00	&      0.00 & 	0.00	& 	1.00	& 	1.00	& 	1.00	& 	1.00	& 	0.02	& 	0.01	& 	0.02	& 	0.02	 \\
\texttt{net12} 			& 1.86e+3	& 2.18e+3	& 2.49e+3	& 2.49e+3	&    19.20	& 	0.00	& 	0.46	& 	0.56	& 	0.32	& 	0.34	&   160.89	&     71.24	&   206.81	&     71.79	 \\
\texttt{p2756} 			& 5.29e+3	&  207.07	& 	inf	& 	inf	&    59.60	&     20.80	& 	0.27	& 	0.38	& 	0.00	& 	0.00	&    23.31	&     11.57	&     20.31	&     12.53	 \\
\texttt{protfold} 			& 	inf	& 	inf	&      inf	&    22.58	&    45.20	&     32.30	& 	0.00	& 	0.00	& 	0.00	& 	0.09	&    33.15	&     19.37	&     27.09	&     11.15	 \\
\texttt{sp97ar} 			&    16.29	&    13.36	&   12.21	&    13.98	&    17.40	&     15.40	& 	0.99	& 	1.00	& 	1.00	& 	1.00	&    30.46	&     35.20	&     35.64	&     39.61	 \\ \hline
\end{tabular}
\end{table*}

\begin{table*}[h]
\caption{Detailed results for different solvers on class B problems. The ``inf'' means that the 
problem has not been solved to feasibility.}
\label{table:classB}
\centering
\begin{tabular}{|l|rrrrrr|rrrr|rrrr|}  \hline
\multirow{2}{*}{{\bf{Class B}}} 	& \multicolumn{6}{c|}{{\bf{gap}}} &  \multicolumn{4}{c|}{{\bf{success}}}	&  \multicolumn{4}{c|}{{\bf{time}}} 	 \\  
					& AFP & FP & AFP-1  & FP-1 & FP2.0 & FPC 
					& AFP & FP & AFP-1 & FP-1 & 
					AFP & FP & AFP-1 & FP-1 \\  \hline			
\texttt{10teams} 		&	0.65	& 	inf	& 	0.65	& 	0.65	& 	4.50	& 	3.20	& 	0.04	& 	0.00	& 	0.18	& 	0.19	&    30.98	& 	9.93	&    11.74	& 	6.78	 \\ 
\texttt{22433} 			& 	0.00	& 	0.00	& 	0.00	& 	0.00	& 	0.00	& 	0.00	& 	0.49	& 	0.54	& 	0.76	& 	0.81	& 	4.22	& 	1.14	& 	0.80	& 	0.37	 \\ 
\texttt{23588} 			& 	2.22	& 	1.59	& 	1.43	& 	0.87	& 	1.10	& 	1.00	& 	0.30	& 	0.51	& 	0.45	& 	0.94	& 	3.68	& 	0.84	& 	1.22	& 	0.31	 \\ 
\texttt{bienst1} 			& 	0.53	& 	0.00	& 	0.00	& 	2.14	& 	0.00	& 	0.00	& 	1.00	& 	1.00	& 	1.00	& 	1.00	& 	1.91	& 	0.46	& 	0.91	& 	0.48	 \\ 
\texttt{bienst2} 			& 	6.83	& 	0.73	& 	3.48	& 	3.33	& 	9.30	& 	2.10	& 	1.00	& 	1.00	& 	1.00	& 	1.00	& 	2.06	& 	0.50	& 	0.93	& 	0.52	 \\ 
\texttt{dano3mip} 		& 	1.91	& 	1.93	& 	2.38	& 	1.48	& 	6.90	& 	6.10	& 	1.00	& 	1.00	& 	1.00	& 	1.00	&    27.48	&    33.09	&    26.70	&     22.34	 \\ 
\texttt{markshare1} 		& 1.32e+4	& 1.21e+4	& 6.80e+3	& 1.40e+4	& 2.19e+4	& 1.22e+4	& 	1.00	& 	1.00	& 	1.00	& 	1.00	& 	0.50	& 	0.10	& 	0.25	& 	0.08	 \\ 
\texttt{markshare2} 		& 1.01e+4	& 1.67e+4	& 1.18e+4	& 1.51e+4	& 2.67e+4	& 2.56e+4	& 	1.00	& 	1.00	& 	1.00	& 	1.00	& 	0.48	& 	0.11	& 	0.27	& 	0.10	 \\ 
\texttt{mcf2} 			& 	7.61	&	1.26	& 	4.77	& 	1.26	& 	1.10	& 	9.60	& 	0.13	& 	0.60	& 	0.26	& 	0.56	& 	8.36	& 	1.73	& 	4.27	&      	1.79	 \\
\texttt{neos-1053234} 	& 	0.00	&    76.00	& 	5.00	&    83.00	&    68.20	&    27.60	& 	1.00	& 	1.00	& 	1.00	& 	1.00	&    13.46	& 	5.92	& 	8.40	& 	5.64	 \\ 
\texttt{neos-1120495} 	&  	0.90	& 	2.40	& 	2.99	& 	2.99	&    21.30	&    18.60	& 	0.60	& 	1.00	& 	0.62	& 	0.99	&    33.37	& 	7.75	& 	8.48	&    	5.86	 \\ 
\texttt{neos-1121679} 	&  731.25	&  731.25	&  331.25	&  781.25	& 1.51e+3	&  700.00	& 	1.00	& 	1.00	& 	1.00	& 	1.00	& 	0.43	& 	0.09	& 	0.24	&     0.084	 \\ 
\texttt{neos-522351} 		&	0.60	& 	2.00	& 	0.51	& 	1.38	&    19.20	& 	0.00	& 	1.00	& 	1.00	& 	1.00	& 	1.00	& 	9.87	& 	2.30	& 	4.60	& 	2.27	 \\  
\texttt{neos-582605} 		&  	0.00	& 	0.00	& 	0.00	& 	0.00	& 	0.00	& 	0.00	& 	0.30	& 	0.29	& 	0.27	& 	0.28	&    13.45	& 	5.59	& 	6.96	& 	5.18	 \\ 
\texttt{neos-583731} 		&   	inf	& 	inf	& 	inf	& 	inf	& 	0.00	& 	0.00	& 	0.00	& 	0.00	& 	0.00	& 	0.00	&    15.28	& 	5.63	& 	5.25	&   	2.46	 \\ 
\texttt{neos-631164} 		& 	9.43	& 	8.91	& 	9.91	& 	9.74	& 	7.00	& 	8.90	& 	1.00	& 	1.00	& 	1.00	& 	1.00	& 	1.52	& 	0.44	& 	0.65	&   	0.41	 \\ 
\texttt{neos-631517} 		& 	4.12	& 	4.01	& 	4.34	& 	4.03	& 	7.50	& 	8.70	& 	1.00	& 	1.00	& 	1.00	& 	1.00	& 	1.36	& 	0.37	& 	0.55	& 	0.32	 \\ 
\texttt{neos-777800} 		& 	0.00	& 	0.00	& 	0.00	& 	0.00	& 	0.00	& 	0.00	& 	0.00	& 	1.00	& 	0.88	& 	1.00	&    21.99	& 	5.19	& 	8.23	& 	5.97	 \\ 
\texttt{neos-831188} 		& 	4.60	& 	inf	&    10.34	& 	inf	& 	1.80	& 	1.60	& 	0.07	& 	0.00	& 	0.14	& 	0.00	&    76.64	&    29.79	&    42.87	&     26.35	 \\ 
\texttt{neos-839838} 		& 	1.71	&    29.16	& 	1.27	&    21.84	&    12.90	& 	9.30	& 	1.00	& 	1.00	& 	1.00	& 	1.00	&    41.63	&    31.51	&    36.64	&     33.23	 \\ 
\texttt{neos-839859} 		&    34.80	&    76.19	&    45.91	&    64.11	&    61.80	&    13.00	& 	1.00	& 	1.00	& 	1.00	& 	1.00	&    10.47	& 	3.17	& 	4.28	& 	2.95	 \\ 
\texttt{neos-886822} 		& 	5.41	& 	2.29	& 	5.75	& 	1.51	&    24.20	&    22.40	& 	1.00	& 	1.00	& 	1.00	& 	1.00	& 	5.94	& 	2.06	& 	2.72	& 	2.00	 \\ 
\texttt{neos-911880} 		&    43.02	&    48.57	&    51.07	&    37.52	&  152.30	&    11.20	& 	1.00	& 	1.00	& 	1.00	& 	1.00	& 	2.35	& 	0.64	& 	1.00	&     	0.59	 \\ 
\texttt{neos-911970} 		&    60.61	&    65.34	&    20.27	&    85.77	&  115.00	&    11.40	& 	1.00	& 	1.00	& 	1.00	& 	1.00	& 	2.52	& 	0.65	& 	1.02	& 	0.62	 \\  
\texttt{neos-912015} 		& 	inf	& 	inf	&    21.43	& 	inf	&    28.60	&    21.40	& 	0.00	& 	0.00	& 	0.06	& 	0.00	&     18.31	& 	5.67	& 	8.21	& 	2.67	 \\ 
\texttt{opt1217} 	     		& 	0.00	& 	0.00	& 	0.00	& 	0.00	& 	0.00	& 	0.00	& 	1.00	& 	1.00	& 	1.00	& 	1.00	& 	0.77	& 	0.19	& 	0.32	& 	0.18	 \\ 
\texttt{pk1} 		    	&    54.55	&  127.27	&    72.73	&  100.00	&  272.70	&    63.60	& 	1.00	& 	1.00	& 	1.00	& 	1.00	& 	0.48	& 	0.12	& 	0.24	& 	0.10	 \\ 
\texttt{swath} 	    		&    32.89	& 1.49e+3	& 	inf	& 	inf	&    65.20	&    35.60	& 	0.31	& 	0.10	& 	0.00	& 	0.00	&    70.68	&    52.81	&    44.14	&     47.97	 \\ 
\texttt{tr12-30} 	    		&    62.69	&    20.39	&    62.44	&    23.28	&    23.50	&    22.60	& 	0.64	& 	1.00	& 	0.62	& 	1.00	& 	5.70	& 	0.77	& 	2.25	&     0.72	 \\ 
\texttt{vpm2} 	    		&    16.36	&    10.91	&    47.27	& 	9.09	&    12.70	& 	3.60	& 	0.44	& 	1.00	& 	0.18	& 	1.00	& 	3.11	& 	0.19	& 	1.60	& 	0.17	 \\  \hline
\end{tabular}
\end{table*}

\begin{table*}[h]
\caption{Detailed results of different solvers on class C problems. The ``inf'' means that the 
problem has not been solved to feasibility.}
\label{table:classC}
\centering
\begin{tabular}{|l|rrrrrr|rrrr|rrrr|}  \hline
\multirow{2}{*}{{\bf{Class C}}} 	& \multicolumn{6}{c|}{{\bf{gap}}} &  \multicolumn{4}{c|}{{\bf{success}}}	&  \multicolumn{4}{c|}{{\bf{time}}} 	 \\  
					& AFP & FP & AFP-1  & FP-1 & FP2.0 & FPC 
					& AFP & FP & AFP-1 & FP-1 & 
					AFP & FP & AFP-1 & FP-1 \\  \hline
\texttt{a1c1s1}			&   88.48	& 	inf	&   88.48	& 	inf	&    38.60	&   21.80	& 	1.00	& 	0.00	& 	1.00	& 	0.00	& 	2.94	&    15.70	& 	2.58	& 	15.66 \\  
\texttt{aligninq} 			&	inf	& 	1.14	& 	1.11	& 	1.33	& 	0.90	&    1.30	& 	0.00	& 	0.14	& 	0.09	& 	0.28	&    10.81	& 	7.90	& 	9.18	& 	6.22	 \\ 
\texttt{danoint}			&   11.16	& 	2.53	& 	5.73	& 	1.26	& 	5.10	&   10.40	& 	0.20	& 	0.71	& 	0.31	& 	0.66	& 	4.05	& 	1.82	& 	4.07	& 	1.55	 \\  
\texttt{liu} 				&   104.86	&   109.01	&    90.27	&    115.5	&   169.60	&   93.00	& 	0.90	& 	0.77	& 	0.94	& 	0.96	& 	5.20	& 	3.98	& 	4.76	& 	2.70	 \\ 
\texttt{mas74} 			&    12.51	& 	4.67	& 	8.23	&    11.46	&    43.10	&    9.10	& 	1.00	& 	1.00	& 	1.00	& 	1.00	& 	0.11	& 	0.06	& 	0.11	& 	0.06	 \\ 
\texttt{mas76} 			&	1.39	& 	0.00	& 	1.90	& 	0.00	& 	6.30	&    2.30	& 	1.00	& 	1.00	& 	1.00	& 	1.00	& 	0.06	& 	0.04	& 	0.06	& 	0.03	 \\ 
\texttt{modglob}			&	1.45	&     11.59	& 	1.93	&    10.63	& 6.50e+6	& 6.40e+6	& 	1.00	& 	1.00	& 	1.00	& 	1.00	& 	0.09	& 	0.17	& 	0.08	& 	0.16	 \\  
\texttt{neos-1053591}	&	3.02	& 	2.30	& 	2.50	& 	inf	& 	2.00	& 	0.10	& 	0.31	& 	0.10	& 	0.09	& 	0.00	& 	6.63	& 	4.64	& 	6.28	& 	4.13	 \\  
\texttt{neos-1112782} 	&  228.67	& 1.96e+3	&  233.92	& 1.92e+3	& 3.76e+3	&   52.40	& 	1.00	& 	1.00	& 	1.00	& 	1.00	& 	7.96	&   12.15  & 	7.37	& 	5.64	 \\ 
\texttt{neos-1112787}	&  232.74	& 2.91e+3	& 2.54e+02& 3e+03	& 3.37e+3	&   37.70	& 	1.00	& 	1.00	& 	1.00	& 	0.99	& 	5.72	& 	3.92	& 	10.7	& 	7.03	 \\  
\texttt{neos-1173026} 	&   31.00	&    32.00	&   31.00	&    37.00	&   35.50	&   20.60	& 	1.00	& 	1.00	& 	1.00	& 	1.00	& 	3.03	& 	1.66	& 	2.96	& 	1.52	 \\ 
\texttt{neos-430149} 		&  613.43	&   152.83	& 	inf	& 	inf	&   485.40	&  164.70	& 	0.11	& 	0.21	& 	0.00	& 	0.00	& 	4.32	& 	2.44	& 	3.50	& 	2.14	 \\ 
\texttt{neos-598183} 		&  129.97	&    15.08	& 	inf	& 	1.22	&   103.40	&    0.60	& 	0.21	& 	0.10	& 	0.00	& 	0.11	& 	8.47	& 	6.90	& 	9.19	& 	6.04	 \\ 
\texttt{neos-603073} 		&	0.24	& 	0.01	& 	0.22	& 	0.02	&   31.30	&    0.40	& 	0.90	& 	0.60	& 	0.95	& 	0.60	& 	3.82	& 	3.60	& 	3.37	& 	3.47	 \\ 
\texttt{neos-631694} 		&	inf	& 	inf	& 	inf	& 	inf	& 	inf	&    0.00	& 	0.00	& 	0.00	& 	0.00	& 	0.00	&    22.66	&   23.08	&    19.64	&     19.71	 \\ 
\texttt{neos-803219}		&	2.71	& 	4.85	& 	inf	& 	inf	& 	7.80	&    4.50	& 	0.24	& 	0.23	& 	0.00	& 	0.00	& 	4.10	& 	2.44	& 	3.41	& 	2.11	 \\  
\texttt{neos-803220} 		&	2.66	& 	2.87	& 	inf	&   12.33	& 	4.70	&    4.40	& 	0.42	& 	0.64	& 	0.00	& 	0.47	& 	3.07	& 	2.17	& 	3.41	& 	1.98	 \\ 
\texttt{neos-806323} 		&	inf	& 	inf	& 	inf	& 	inf	&    44.70	&    8.70	& 	0.00	& 	0.00	& 	0.00	& 	0.00	& 	7.27	& 	4.55	& 	5.51	& 	3.45	 \\ 
\texttt{neos-807639} 		&	1.34	& 	4.03	& 	1.34	& 	5.37	& 	2.70	&    6.70	& 	0.33	& 	0.19	& 	0.14	& 	0.06	& 	5.21	& 	3.80	& 	5.84	& 	3.60	 \\ 
\texttt{neos-807705}		&    10.43	&    11.95	& 	inf	& 	inf	& 	8.80	&   15.70	& 	0.22	& 	0.07	& 	0.00	& 	0.00	& 	6.58	& 	4.12	& 	5.73	& 	3.53	 \\  
\texttt{neos-848589} 		&    95.77	& 3.79e+4	& 5.03e+4& 265.85& 6.79e+4	&  526.80	&      0.92	&      1.00	&      0.92	& 	1.00	&3.14e+3	& 2.3e+3  & 3.9e+3	& 2.7e+3	 \\ 
\texttt{neos-863472}		&	0.00	& 	0.00	& 	inf	& 	inf	&    11.10	&    0.70	& 	0.46	& 	0.36	& 	0.00	& 	0.00	& 	2.43	& 	1.78	& 	2.92	& 	1.91	 \\  
\texttt{neos-880324}		&	inf	&     inf	& 	inf	& 	inf	& 	inf	&    0.00	& 	0.00	& 	0.00	& 	0.00	& 	0.00	& 	1.81	& 	1.01	& 	1.29	& 	0.84	 \\  
\texttt{neos-912023} 		&	inf	&     inf	& 	inf	&   38.46	&   15.40	&    0.00	& 	0.00	& 	0.00	& 	0.00	& 	0.04	& 	8.96	& 	5.48	& 	7.11	& 	2.41	 \\ 
\texttt{neos-913984} 		&	inf	&     inf	&   inf	&   inf	&    0.00	&    0.00	&      0.00	& 	0.00	& 	0.00	& 	0.00	& 4.6e+3	& 10e+3	& 5.9e+3	&   11.6e+3\\ 
\texttt{neos-916792} 		&    2.36	&    3.90	&   10.05	&   -2.48	&   10.80	&    7.20	& 	0.43	& 	0.30	& 	0.54	& 	0.34	&    22.00	&    17.97	&    19.24	&    13.17	 \\ 
\texttt{neos17} 			&    43.00	&   38.00	&   43.00	&   37.00	&   51.60	&   23.00	& 	1.00	& 	1.00	& 	1.00	& 	1.00	& 	2.25	& 	1.26	& 	2.09	& 	1.22	 \\ 
\texttt{neos2} 			&	inf	& 	inf	& 	inf	& 	inf	& 	inf	&  197.70	& 	0.00	& 	0.00	& 	0.00	& 	0.00	&    14.71	& 	8.58	& 	9.91	& 	7.96	 \\ 
\texttt{neos4} 			&	inf	& 	inf	&   inf	&      inf	& 	inf	&    5.10	& 	0.00	& 	0.00	&.      0.00	& 	0.00	& 1.5e+2	& 2e+2     & 9.3e+2	& 11.2e+2	 \\ 
\texttt{set1ch}			&   75.73	&   20.39	&   82.31	&   20.24	& 47.70	&    9.00	& 	1.00	& 	1.00	& 	1.00	& 	1.00	& 	0.21	& 	0.40	& 	0.19	& 	0.36	 \\  \hline 
\end{tabular}
\end{table*}

\end{document}